\documentclass{amsart}%
\usepackage{amsmath}
\usepackage{amsfonts}
\usepackage{amssymb}
\usepackage{graphicx}%
\usepackage{xcolor}
\usepackage{hyperref}
\newtheorem{lemma}{Lemma}[section]
\newtheorem{theorem}[lemma]{Theorem}
\newtheorem{remark}[lemma]{Remark}

\newtheorem{coro}[lemma]{Corollary}
\newtheorem{definition}[lemma]{Definition}
\newtheorem{example}[lemma]{Example}

\title[Almost Periodic Solutions of Lattice Dynamical Systems \dots ]{Almost Periodic Solutions of Lattice Dynamical Systems with Monotone Nonlinearity}

\author{David Cheban}
\address[D. Cheban]{State University of Moldova\\
Vladimir Andrunachievici Instiy=tute of Mathematics and Computer Science\\
Laboratory of Differential Equations\\str. Academiei 5\\
MD--2028 Chi\c{s}in\u{a}u, Moldova
} \email[D.
Cheban]{david.ceban@usm.md, davidcheban@yahoo.com}
\author{Andrei Sultan}
\address[A. Sultan]{%
State University of Moldova\\Vladimir Andrunachievici Institute of
Mathematics and Computer Science\\ Laboratory of Differential Equations\\str. Academiei 5\\
MD--2028 Chi\c{s}in\u{a}u, Moldova} \email[A.
Sultan]{andrew15sultan@gmail.com}

\date{\today}

\begin{document}

\begin{abstract}
The aim of this paper is studying the problem of almost
periodicity of almost periodic lattice dynamical systems of the
form $u_{i}'=\nu (u_{i-1}-2u_i+u_{i+1})-\lambda
u_{i}+F(u_i)+f_{i}(t)\ (i\in \mathbb Z,\ \lambda >0)$. We prove
the existence a unique almost periodic solution of this system if
the nonlinearity $F$ is monotone.
\end{abstract}

\maketitle

\section{Introduction}\label{Sec1}

Denote by $\mathbb R :=(-\infty,\infty)$, $\mathbb Z :=\{0,\pm
1,\pm 2,\ldots\}$ and $\ell_{2}$ the Hilbert space of all
two-sided sequences $\xi =(\xi_{i})_{i\in \mathbb Z}$ ($\xi_{i}\in
\mathbb R$) with
\begin{equation}\label{eqI_01}
\sum\limits_{i\in \mathbb Z}|\xi_{i}|^{2}<+\infty \nonumber
\end{equation}
and equipped with the scalar product
\begin{equation}\label{eqI_2}
\langle \xi,\eta\rangle :=\sum\limits_{i\in \mathbb
Z}\xi_{i}\eta_{i} .\nonumber
\end{equation}

Let $(X,\rho)$ be a complete metric space with the distance
$\rho$, $C(\mathbb R,X)$ be the space of all continuous functions
$f:\mathbb R\to X$ equipped with the distance
\begin{equation}\label{eqI_3}
d(f_1,f_2):=\sup\limits_{L>0}\min\{\max\limits_{|t|\le
L}\rho(f_{1}(t),f_{2}(t)),L^{-1}\}.
\end{equation}
The metric space $(C(\mathbb R,X),d)$ is complete and the distance
$d$, defined by (\ref{eqI_3}), generates on the space $C(\mathbb
R,X)$ the compact-open topology.

Let $h\in \mathbb R$, $f\in C(\mathbb R,X)$, $f^{h}(t):=f(t+h)$
for all $t\in \mathbb R$ and $\sigma :\mathbb R\times C(\mathbb
R,X)\to C(\mathbb R,X)$ be a mapping defined by
$\sigma(h,f):=f^{h}$ for all $(h,f)\in \mathbb R\times C(\mathbb
R,X)$. Then \cite[Ch.I]{Che_2015} the triplet $(C(\mathbb
R,X),\mathbb R,\sigma)$ is a shift dynamical system (or Bebutov's
dynamical system) on he space $C(\mathbb R,X)$. By $H(f)$ the
closure in the space $C(\mathbb R,X)$ of $\{f^{h}|\ h\in \mathbb
R\}$ is denoted.

\begin{definition}\label{defAPF02} A function $f \in C(\mathbb
T,X)$ is said to be Lagrange stable if the motion $\sigma(t,f)$ is
so in the shift dynamical system $(C(\mathbb T,X),\mathbb
T,\sigma)$, i.e., $H(f)$ is a compact subset of $C(\mathbb R,X)$.
\end{definition}

\begin{lemma}\label{lAPF02}  {\rm(\cite[Ch.IV,p.236]{Bro79},\cite[Ch.III]{Sel_1971},\cite[Ch.IV]{sib})}
A function $\varphi \in C(\mathbb T,X)$ is Lagrange stable if and
only if the following conditions are fulfilled:
\begin{enumerate}
\item the set $\varphi(\mathbb T):=\{\varphi(t)|\ t\in \mathbb
T\}$ is precompact in $X$; \item the function $\varphi$ is
uniformly continuous on $\mathbb T$.
\end{enumerate}
\end{lemma}

Recall that a subset $A\subset \mathbb R$ is called relatively
dense if there exits a positive number $l$ such that
\begin{equation}\label{eqRD1}
A\bigcap [a,a+l]\not= \emptyset \nonumber
\end{equation}
for all $a\in \mathbb R$.

A function $f\in C(\mathbb R,X)$ is called almost periodic
\cite{Che_2020,Lev-Zhi}, if for every positive number $\varepsilon$
the set
\begin{equation}\label{eqRD2}
\mathcal F(f,\varepsilon):=\{\tau \in \mathbb R |\
\rho(f(t+\tau),f(t))<\varepsilon\ \ \mbox{for all}\ t\in \mathbb
R\} \nonumber
\end{equation}
is relatively dense.

In this paper we study the problem of existence at least one
almost periodic solution of the systems
\begin{equation}\label{eqI1}
u_{i}'=\nu (u_{i-1}-2u_i+u_{i+1})-\lambda u_{i}+F(u_i)+f_{i}(t)\
(i\in \mathbb Z),
\end{equation}
where $\lambda >0$, $F\in C(\mathbb R, \mathbb R)$ and $f\in
C(\mathbb R,\ell_{2})$ ($f(t):=(f_{i}(t))_{i\in \mathbb Z}$ for
all $t\in \mathbb R$) is an almost periodic function.

The system (\ref{eqI1}) can be considered as a discrete (see, for
example, \cite{BLW_2001,HK_2023} and the bibliography therein)
analogue of a reaction-diffusion equation in $\mathbb R$:
\begin{equation}\label{eqI1.1}
 \frac{\partial{u}}{\partial{t}} = D\frac{\partial^{2}{u}}{\partial^{2}{x}}-\lambda u + F(u) +
 f(t,x),\nonumber
\end{equation}
where grid points are spaced $h$ distance apart and $\nu =
D/h^{2}$.

This study continues the authors work \cite{CS_2025} devoted to
the study the problem of existence of compact global attractor for
(\ref{eqI1}) and the work \cite{CS} dedicated to the study the
invariant sections of monotone nonautonomous dynamical systems and
their applications to differen classes of evolution equations. The
invariant sections play a very important role in the study the
problem of existence of almost periodic (respectively, almost
automorphic, recurrent and Poisson stable) solutions of
differential equations.

The paper is organized as follows. In the second section we
collect some notions and facts from dynamical systems (both
autonomous and nonautonomous) which we use in this paper. The
third second section is dedicated to the proof that under some
conditions the equation (\ref{eqI1}) generates a cocycle which
plays a very important role in the study of the asymptotic
properties of the equation (\ref{eqI1}). In the fourth section we
show that under some conditions there exists a compact global
attractor for the equation (\ref{eqI1}). The fifth section is
dedicated to the study the invariant sections of the cocycle
generated by the equation (\ref{eqI1}). In the sixth section we
study the structure of the compact global attractor for the
equation (\ref{eqI1}). Namely, we show that the equation
(\ref{eqI1}) is convergent, i.e., it admits a compact global
attractor $I =\{I_{y}|\ y\in Y\}$ such that every set $I_{y}$
consists of a single point. The seventh section is dedicated to
the application of our general results from Sections
\ref{Sec3}-\ref{Sec6}. Namely we are analyzing an example of
equation of the form (\ref{eqI1}) which illustrate our general
results. Finally, in Section \ref{Sec8} we give some
generalization of our results (Theorem \ref{thAPM1} and Corollary
\ref{corAPM1}) for almost automorphic and recurrent lattice
dynamical systems,

\section{Preliminary}\label{Sec2}

Below we give some notions, notation and facts from the theory of
dynamical systems \cite{Bro79,Che_2020,Che_2024,Sch72,Sch85,sib}
which we will use in this paper.

Let $(X,\rho_{X})$ and $(Y,\rho_{Y})$ be two complete metric
spaces with the distance $\rho_{X}$ and $\rho_{Y}$
respectively\footnote{In what follows, in the notation $\rho_{X}$
(respectively, $\rho_{Y}$), we will omit the index $X$
(respectively, $Y$) if this does not lead to a misunderstanding.},
$\mathbb R:=(-\infty,+\infty)$ and $\mathbb Z :=\{0,\pm 1, \pm 2,
\ldots \}$, $\mathbb T =\mathbb R$ or $\mathbb R_{+}$. Let
$(X,\mathbb R_{+},\pi)$ (respectively, $(Y,\mathbb R, \sigma )$)
be an autonomous one-sided (respectively, two-sided) dynamical
system on $X$ (respectively, on $Y$).

\begin{definition}\label{def1.0.1}
A triplet $(X,\mathbb{T},\pi)$, where $\pi:\mathbb{T}\times X\to
X$ is a continuous mapping satisfying the  conditions $
\pi(0,x)=x$ and $\pi(s,\pi(t,x))=\pi(s+t,x)$ (for every $x\in X$ and
$t,s\in\mathbb T$) is called a dynamical system.
\end{definition}

If $\mathbb{T}=\mathbb{R}$ (respectively, $\mathbb{R_{+}})$), then
$(X,\mathbb{T},\pi)$ is called a group (respectively, semigroup)
dynamical system.

\begin{definition}\label{def1.0.2}
The function $\pi(\cdot,x):\mathbb{T} \to X$ is called a motion
passing through the point $x$ at the initial moment $t=0$ and the
set $\Sigma_x:= \pi(\mathbb{T},x)$ is called a trajectory of this
motion.
\end{definition}

\begin{definition}\label{def1.0.7}
A point $x\in X$ is called a stationary (respectively,
\emph{$\tau$-periodic}, $\tau
>0,\ \tau \in \mathbb T$) \emph{point} if $xt=x$ (respectively, $x\tau = x$)
for all $t\in \mathbb T$, where $xt:=\pi (t,x)$.
\end{definition}

An $m$-dimensional torus is denoted by $\mathcal T^m:=\mathbb
R^m/2{\pi}\mathbb Z^{m}.$ Let $(\mathcal T^m,\mathbb T,\sigma)$ be
an irrational winding of $\mathcal T^m$ with the frequency $\nu
=(\nu_1,\nu_2,\ldots,\nu_m)\in\mathbb R^{m}$, i.e.,
$\sigma(t,v):=(v_1+\nu_1 t (mod\ 2\pi),v_2+\nu_2 t(mod\ 2
\pi),\ldots,v_m+\nu_m t(mod\ 2\pi))$ for all $t\in \mathbb T$ and
$v=(v_1,v_2,\ldots,v_m )\in \mathcal T ^m$, where the numbers
$\nu_1,\nu_2,\ldots,\nu_m$ are rational independent.

\begin{definition}\label{defQP}
A point $x\in X$ (respectively, a motion $\pi(t,x)$) is called
quasi periodic with the frequency $\mu
:=(\mu_1,\mu_2,\ldots,\mu_m)\in \mathbb R^m$ if there exists a
continuous function $\mathfrak F :\mathcal T^m \to X$ such that
$\mathfrak F(v_0)=x$ for some $v_0\in \mathcal T^{m}$ and
$\mathfrak F (\sigma(t,v))=\pi(t,\mathfrak F(v))$ for every
$(t,v)\in \mathbb T \times \mathcal T^{m},$ where
 $(\mathcal T^m,\mathbb
T,\sigma)$ is an irrational winding of torus $\mathcal T^m$ with
the frequency $\mu =(\mu_1,\mu_2,\ldots,\mu_m)$.
\end{definition}

\begin{definition}\label{defAPM1} A point $x\in X$ (respectively,
the motion $\pi(t,x)$) is said to be almost periodic
\cite[Ch.I]{Che_2020} if for arbitrary positive number
$\varepsilon$ the set
\begin{equation}\label{eqAPM1}
\mathcal F(x,\varepsilon):=\{\tau \in \mathbb R|\
\rho(\pi(t+\tau,x),\pi(t,x))< \varepsilon\ \forall \ t\in \mathbb
R\}\nonumber
\end{equation}
is relatively dense.
\end{definition}

\begin{remark}\label{remAP2} Every quasi-periodic point is almost
periodic \cite[Ch.I]{Che_2024}.
\end{remark}

\begin{lemma}\label{lAP01} \cite{Che_2024, Sch72} Assume that $f\in C(\mathbb R,X)$ is a
continuous function $f:\mathbb R\to X$ and $((C(\mathbb
R,X),\mathbb R,\sigma)$ is the shift dynamical system on the space
$C(\mathbb R,X)$.

The following statements are equivalent:
\begin{enumerate}
\item the function $f$ is stationary (respectively,
$\tau$-periodic, quasi-periodic or almost periodic); \item the
motion $\sigma(t,f)$ generated by the function $f$ in the shift
dynamical system $(C(\mathbb R,X),\mathbb R,\sigma)$ is stationary
(respectively, $\tau$-periodic, quasi-periodic or almost
periodic).
\end{enumerate}
\end{lemma}

\begin{definition}\label{def1.0.17}
Let $(X,\mathbb{T}_1,\pi)$ and $(Y,\mathbb{T}_2,\sigma)$ \
($\mathbb T_{1},\mathbb T_{2}\in \{\mathbb R_{+},\mathbb R\},\
\mathbb{R_{+}}\subseteq\mathbb{T}_1\subseteq\mathbb{T}_2
\subseteq\mathbb{R}$) be two dynamical systems. A mapping $h:X\to
Y$ is called a homomorphism of dynamical system
$(X,\mathbb{T}_1,\pi)$ into $(Y,\mathbb{T}_2,\sigma)$ if the
mapping $h$ is continuous and $h(\pi(x,t))=\sigma(h(x),t)$ (for
every $t\in\mathbb{T}_1$ and $x\in X$).
\end{definition}

\begin{definition}\label{def1.0.18} Let $(X,h,Y)$ be a bundle \cite[Ch.I]{hew}.
The triplet $\langle (X,$ $\mathbb{T}_1,$ $\pi),$
$\,(Y,\mathbb{T}_2,\sigma), $ $\,h\rangle $, where $h$ is a
homomorphism from $(X,$ $\mathbb{T}_1,$ $\pi)$ on $(Y,$
$\mathbb{T}_2,$ $\sigma)$ is called a nonautonomous dynamical
system.
\end{definition}

\begin{definition}\label{def1.0.19}
The triplet $\langle W, \varphi, (Y,\mathbb{T}_2,\sigma)\rangle
$\index{$\langle W, \varphi, (Y,\mathbb{T}_2,\sigma)\rangle $} (or
shortly $\varphi$), where $(Y,\mathbb{T}_2,\sigma)$ is a dynamical
system on $Y$,  $W$ is a complete metric space and $\varphi$ is a
continuous mapping from $\mathbb{T}_1\times W\times Y$ into $W$,
possessing the following properties:
\begin{enumerate}
\item[a.] $\varphi(0,u,y)=u$ (for all $u\in W, y\in Y)$; \item[b.]
$\varphi(t+\tau,u,y)= \varphi(\tau,\varphi(t,u,y),\sigma(t,y))$
for all $t,\tau\in\mathbb{T}_1,\, u\in W$ and $y\in Y$,
\end{enumerate}
is called \cite{Arn,Sel_1971} a \emph{cocycle} over
$(Y,\mathbb{T}_2,\sigma)$ with the fibre $W$.
\end{definition}

\begin{definition}\label{def1.0.20}
Let $X:= W\times Y$ and we define a mapping $\pi: X\times
\mathbb{T}_1\to X$ as following:
$\pi((u,y),t):=(\varphi(t,u,y),\sigma(t,y))$ for every $(u,y)\in
X$ and $t\in \mathbb T_1$ (i.e., $\pi=(\varphi,\sigma)$). Then it
easy to see that $(X,\mathbb{T}_1,\pi)$ is a dynamical system on
$X$, associated by the cocycle $\varphi$, which is called a
skew-product dynamical system \cite{alek,Sel_1971} and
$h=pr_2:X\to Y$ is a homomorphism from $(X,\mathbb{T}_1,\pi)$ on
$(Y,\mathbb{T}_2,\sigma)$ and, consequently, $\langle
(X,\mathbb{T}_1,\pi),\, (Y,\mathbb{T}_2,\sigma), h\rangle $ is
\emph{a nonautonomous dynamical system, associated/generated by
the cocycle $\varphi$}.
\end{definition}

Thus if we have a cocycle $\langle W, \varphi, (Y,\mathbb{T}_2,
\sigma)\rangle $ over dynamical system $(Y,$ $\mathbb{T}_2,$
$\sigma)$ with the fibre $W$, then it generates a nonautonomous
dynamical system $\langle (X,\mathbb{T}_1,\pi),$\
$(Y,\mathbb{T}_2,\sigma), h\rangle $ ($X:= W\times Y$).

\begin{example}\label{ex1.1.13}
{\rm Let $\mathfrak B$ be a real or complex Banach space with the
norm $|\cdot|$. Let us consider a differential equation
\begin{equation}
u'=f(t,u),\label{eq1.0.6}
\end{equation}
where  $f\in C(\mathbb{R}\times \mathfrak B,\mathfrak B)$. Along
with the equation $(\ref{eq1.0.6})$ we consider its
\emph{$H$-class}\index{$H$-class}
\cite{Bro79,Dem67,Lev-Zhi,Sch72,Sch85}, i.e., the family of
equations
\begin{equation}
v'=g(t,v),\label{eq1.0.7}
\end{equation}
where $g\in H(f):=\overline{\{f_{\tau}:\tau\in \mathbb{R}\}}$,
$f_{\tau}(t,u)=f(t+\tau,u)$ for all $(t,u)\in \mathbb{R}\times
\mathfrak B$ and by bar we denote the closure in
$C(\mathbb{R}\times \mathfrak B,\mathfrak B)$. We will suppose
also that the function $f$ is regular \cite[ChIV]{Sel_1971}, i.e.,
for every equation (\ref{eq1.0.7}) the conditions of existence,
uniqueness and extendability on $\mathbb{R}_{+}$ are fulfilled.
Denote by $\varphi(\cdot,v,g)$ the solution of the equation
(\ref{eq1.0.7}) passing through the point $v\in \mathbb B$ at the
initial moment $t=0$. Then from the general properties of
solutions of ordinary differential equations (ODEs) it follows
that the mapping $\varphi:\mathbb{R}_{+}\times \mathfrak B\times
H(f)\to \mathfrak B$ is well defined and it satisfies the
following conditions (see, for example, \cite[ChIV]{Bro79},
\cite{DK_1970} and \cite[ChIV]{Sel_1971}):
\begin{enumerate}
\item[$1)$] $\varphi(0,v,g)=v$ for every $v\in \mathfrak B$ and
$g\in H(f)$;
\item[$2)$]
$\varphi(t,\varphi(\tau,v,g),g_{\tau})=\varphi(t+\tau,v,g)$ for
every $ v\in \mathfrak B$, $g\in H(f)$ and $t,\tau \in
\mathbb{R}_{+}$; \item[$3)$] the mapping
$\varphi:\mathbb{R}_{+}\times \mathfrak B\times H(f)\to \mathfrak
B$ is continuous.
\end{enumerate}

Denote by $Y:=H(f)$ and $(Y,\mathbb{R},\sigma)$ a dynamical system
of translations
 on $Y$, induced by the dynamical system of
translations $(C(\mathbb{R}\times \mathfrak B,\mathfrak
B),\mathbb{R},\sigma)$. The triplet $\langle \mathfrak B,\varphi,
(Y,\mathbb{R},\sigma)\rangle $ is a cocycle over
$(Y,\mathbb{R},\sigma)$ with the fibre $\mathfrak B$. Thus the
equation (\ref{eq1.0.6}) generates a cocycle $\langle \mathfrak
B,\varphi, (Y,\mathbb{R},\sigma)\rangle $ and a nonautonomous
dynamical system $\langle (X,\mathbb{R}_{+},\pi),\,
(Y,\mathbb{R},\sigma), h\rangle $, where $X:= \mathfrak B\times
Y$, $\pi:=(\varphi,\sigma)$ and $h:=pr_2:X\to Y$.}
\end{example}

\section{Cocycles generated by lattice dynamical system (\ref{eqI1}).}\label{Sec3}

Consider a non-autonomous system
\begin{equation}\label{eq2.1}
u_{i}'=\nu (u_{i-1}-2u_i+u_{i+1})-\lambda u_{i}+F(u_i)+f_{i}(t)\
(i\in \mathbb Z) .
\end{equation}

Below we use the following conditions.

\emph{Condition }(\textbf{C1}). \label{C1} The function $f\in
C(\mathbb R,\mathfrak B)$ is almost periodic.

\emph{Condition} (\textbf{C2}). \label{C2} The function $F\in
C(\mathbb R,\mathbb R)$ is Lipschitz continuous on bounded sets
and $F(0)=0$.

Denote by $ \widetilde{F}:\ell_{2}\to \ell_{2}$ the Nemytskii
operator generated by $F$, i.e.,
$\widetilde{F}(\xi)_{i}:=F(\xi_{i})$ for all $i\in \mathbb Z$.

\emph{Condition} (\textbf{C3}). \label{C3} The function $F$ is
monotone, i.e., there exists a number $\alpha \ge 0$ such that
\begin{equation}\label{eqAP_02}
(x_1-x_2)(F(x_1)-F(x_2))\leq -\alpha
|x_1-x_2|^2
\end{equation}
for every $x_1,x_2 \in \mathbb R$.


\begin{lemma}\label{lAP1} The following statements hold:
\begin{enumerate}
\item if the function $f$ satisfies the Conditions (\textbf{C2}),
(\textbf{C3}) and $F(0)=0$, then
\begin{equation}\label{eqAP_03}
F(s)s\le -\alpha s^{2}
\end{equation}
for all $s\in \mathbb R$; \item if the function $F$ satisfies the
Condition (\textbf{C3}), then the Nemytskii operator
$\widetilde{F}$ generated by $F$ possesses the following property
\begin{equation}\label{eqAP04}
\langle u^1-u^2,\widetilde{F}(u^1)-\widetilde{F}(u^2)\rangle \le
-\alpha \|u^1-u^2\|^{2} \nonumber
\end{equation}
for every $u^1,u^2\in \ell_{2}$.
\end{enumerate}
\end{lemma}
\begin{proof}
1. Putting $x_1=s$ and $x_2=0$ in (\ref{eqAP_02}) we obtain
(\ref{eqAP_03}).

2. Let $u^1=(u^{1}_{i})_{i\in \mathbb Z}$ and
$u^2=(u^{2}_{i})_{i\in \mathbb Z}$ be two elements of $\ell_{2}$.
Then, using the monotonicity of $F$, we have
\begin{eqnarray}\label{eqAP5}
& \langle u^1-u^2,\widetilde{F}(u^1)-\widetilde{F}(u^2)\rangle =
\sum\limits_{i\in \mathbb
Z}(u^{1}_{i}-u^{2}_{i})(F(u^{1}_{i})-F(u^{2}_{i}))\le \nonumber
\\
& \sum\limits_{i\in \mathbb Z} -\alpha|u^{1}_{i}-u^{2}_{i}|^{2}=
 -\alpha \|u^1-u^2\|^{2}.\nonumber
\end{eqnarray}
\end{proof}

\begin{definition}\label{defL1.8} A function $F\in C(Y\times \mathfrak B,\mathfrak
B)$ is said to be Lipschitzian on bounded subsets from $\mathfrak
B$ if for every bounded subset $B\subset \mathfrak B$ there exists a
positive constant $L_{B}$ such that
\begin{equation}\label{eqL1.82}
|F(y,v_1)-F(y,v_2)|\le L_{B} |v_1-v_2|
\end{equation}
for all $v_1, v_2 \in B \subset \mathfrak B$.
\end{definition}

\begin{definition}\label{defL2.8}
The smallest constant $L$ (respectively $L_{B}$) with the property
(\ref{eqL1.82}) is called Lipshchitz constant of function $F$
(notation $Lip(F)$ (respectively, $Lip_{B}(F)$)).
\end{definition}

Let $B \subset \mathfrak B$, denote by $CL(Y\times B,\mathfrak B)$
the Banach space of all Lipschitzian functions $F\in C(Y\times
B,\mathfrak B)$ equipped with the norm
\begin{equation}
||F||_{CL}:=\max\limits_{y\in Y}|F(y,0)|+Lip_{B}(F).\nonumber
\end{equation}

\begin{lemma}\label{l2.2} \cite{BLW_2001} Under the Condition (\textbf{C2})
it is well defined the mapping $\widetilde{F}:\ell_{2}\to
\ell_{2}$ and
\begin{equation}\label{eq2.2}
\|\widetilde{F}(\xi)-\widetilde{F}(\eta)\|\le Lip_{B}(F)\|\xi
-\eta \| \nonumber
\end{equation}
for every $\xi,\eta\in \ell_{2}$, where $\|\cdot\|^{2}:=\langle
\cdot,\cdot \rangle$ and $\|\cdot \|$ is the norm on the space
$\ell_{2}$.
\end{lemma}

For every $u = (u_{i})_{i\in \mathbb Z}$, the discrete Laplace
operator $\Lambda$ is defined \cite[Ch.III]{HK_2023} from
$\ell_{2}$ to $\ell_{2}$ component wise by $\Lambda(u)_{i} =
u_{i-1} - 2u_{i} + u_{i+1}$ ($i\in \mathbb Z$). Define the
bounded linear operators $D^{+}$ and $D^{-}$ from $\ell_{2}$ to
$\ell_{2}$ by $(D^{+}u)_{i} = u_{i+1} - u_{i},\ (D^{-}u)_{i} =
u_{i-1} - u_{i}\ (i\in \mathbb Z)$.

Note that $\Lambda = D^{+}D^{-} = D^{-}D^{+}$ and $\langle D^{-}u,
v\rangle = \langle u, D^{+}v\rangle $ for all $u,v\in \ell_{2}$
and, consequently, $\langle \Lambda u,u \rangle = -|D^{+}u|^{2}\le
0$. Since $\Lambda$ is a bounded linear operator acting on the
space $\ell_{2}$, it generates a uniformly continuous semi-group
$\{e^{\Lambda t}\}_{t\ge 0}$ on $\ell_{2}$.

Under the Conditions (\textbf{C1}) and (\textbf{C2}) the system of
differential equations (\ref{eq2.1}) can be written in the form of
an ordinary differential equation
\begin{equation}\label{eq2.3}
u'=\nu \Lambda u +\Phi (u)+f(t)
\end{equation}
in the Banach space $\mathfrak B=\ell_{2}$, where
$\Phi(u):=-\lambda u +\widetilde{F} (u)$ and
$\Lambda(u)_{i}:=u_{i-1}-2u_{i}+u_{i+1}$ for every $u=(u_i)_{i\in
\mathbb Z}\in \ell_{2}$. Along with the equation (\ref{eq2.3}) we
consider also its $H$-class, i.e., the family of equations
\begin{equation}\label{eq2.3g}
u'=\nu \Lambda u +\Phi (u)+g(t),
\end{equation}
where $g\in H(f)$.

The family of the equations (\ref{eq2.3g}) can be rewritten as
follows
\begin{equation}\label{eq2.3H}
u'=F(\sigma(t,g),u)\ \ (g\in H(f)),\nonumber
\end{equation}
where $F:H(f)\times \ell_{2}\to \ell_{2}$ is defined by
$F(g,u):=\nu \Lambda u+\Phi (u) +g(0)$. It easy to see that
$F(\sigma(t,g),u)=\nu \Lambda u+\Phi (u)+g(t)$ for all $(t,u,g)\in
\mathbb R\times \mathfrak B \times H(f)$.

Let $(Y,\mathbb R,\sigma)$ be a dynamical system on the metric
space $Y$.

\begin{theorem}\label{th_AP1} \cite{CS_2025} Under the Conditions (\textbf{C1})-(\textbf{C3}) the following statements hold:
\begin{enumerate}
    \item for every $(v,g)\in \ell_{2}\times H(f)$ there exists a unique
    solution $\varphi(t,v,g)$ of the equation (\ref{eq2.3g}) passing
    through the point $v$ at the initial moment $t=0$ and defined on
    the semi-axis $\mathbb R_{+}:=[0,+\infty)$; \item
    $\varphi(0,v,g)=v$ for all $(v,g)\in \ell_{2}\times H(f)$; \item
    $\varphi(t+\tau,v,g)=\varphi(t,\varphi(\tau,v,g),g^{\tau})$ for
    every $t,\tau\in \mathbb R_{+}$, $v\in \ell_{2}$ and $g\in H(f)$;
     \item the mapping
    $\varphi :\mathbb R_{+}\times \ell_{2}\times H(f)\to \ell_{2}$
    ($(t,v,g)\to \varphi(t,v,g)$ for all $(t,v,g)\in \mathbb
    R_{+}\times \ell_{2}\times H(f)$) is continuous.
\end{enumerate}
\end{theorem}

\begin{coro}\label{corH1}
Under the conditions of Theorem \ref{th_AP1} the equation
(\ref{eq2.3}) (respectively, the family of equations
(\ref{eq2.3g})) generates a cocycle $\langle
\ell_{2},\varphi,(H(f),\mathbb R,\sigma)\rangle$ over the shift
dynamical system $(H(f),\mathbb R,\sigma)$ with the fibre
$\ell_{2}$.
\end{coro}

\begin{theorem}\label{thAP2} Under the Conditions (\textbf{C1})-(\textbf{C3}) the
cocycle $\langle \ell_{2},\varphi, (H(f),\mathbb R,\sigma)\rangle$
generated by the equation (\ref{eq2.3}) possesses the following
property:
\begin{equation}\label{eqAP_14}
\|\varphi(t,v_2,g)-\varphi(t,v_1,g)\|\le e^{-(\lambda +\alpha)
t}\|v_2-v_1\| \nonumber
\end{equation}
for all $v_1,v_2\in \ell_{2}$, $t\ge 0$ and $g\in H(f)$.
\end{theorem}
\begin{proof}
Let $v_1,v_2\in \ell_{2}$, $g\in H(f)$, and consider the solutions
\(\varphi(t,v_1,g)\) and \(\varphi(t,v_2,g)\) of the equation
(\ref{eq2.3g}). Denote by
\(w(t)=\varphi(t,v_2,g)-\varphi(t,v_1,g)\). We aim to show that
\(\|w(t)\|\le e^{-(\lambda +\alpha) t}\|v_2-v_1\|\). Indeed
$$
w' = \varphi(t,v_2,g)' - \varphi(t,v_1,g)'
$$
$$
= \nu \Lambda \varphi(t,v_2,g) + \Phi(\varphi(t,v_2,g)) + g(t) -
(\nu \Lambda \varphi(t,v_1,g) + \Phi(\varphi(t,v_1,g)) + g(t))
$$
$$
= \nu \Lambda w + \Phi(\varphi(t,v_2,g)) - \Phi(\varphi(t,v_1,g))
$$
$$
= \nu \Lambda w - \lambda \varphi(t,v_2,g) + \widetilde{F}(\varphi(t,v_2,g)) + \lambda \varphi(t,v_1,g) - \widetilde{F}(\varphi(t,v_1,g))
$$
$$
= \nu \Lambda w - \lambda w + \widetilde{F}(\varphi(t,v_2,g)) -
\widetilde{F}(\varphi(t,v_1,g))
$$
and
$$
\langle w, w' \rangle = \langle w, \nu \Lambda w - \lambda w + \widetilde{F}(\varphi(t,v_2,g)) - \widetilde{F}(\varphi(t,v_1,g)) \rangle
$$
$$
= \langle w, \nu \Lambda w \rangle - \lambda \langle w, w \rangle
+ \langle w, \widetilde{F}(\varphi(t,v_2,g)) -
\widetilde{F}(\varphi(t,v_1,g)) \rangle .
$$

Evaluate each term:
$$
\langle w,\nu  \Lambda w \rangle = \nu \langle \Lambda w, w \rangle = -\nu |D^{+}w|^{2} \le 0,
$$
$$
-\lambda \langle w, w \rangle = -\lambda \|w\|^{2},
$$
$$
\langle \varphi(t,v_2,g) - \varphi(t,v_1,g),
\widetilde{F}(\varphi(t,v_2,g)) - \widetilde{F}(\varphi(t,v_1,g))
\rangle \le -\alpha \|w(t)\|^{2} \quad \mbox{(by Lemma
\ref{lAP1})},
$$
where $w(t)=\varphi(t,v_2,g) - \varphi(t,v_1,g)$ for all $t\in
\mathbb R_{+}$. Combining these results, we have:
    \begin{align*}
        \frac{d}{dt}\|w(t)\|^{2} =2\langle w(t), w'(t) \rangle \le -2(\lambda + \alpha) \|w(t)\|^{2}.
    \end{align*}
    By Gronwall's inequality, it follows that
    \begin{align*}
        \|w(t)\|^{2} \le e^{-2(\lambda + \alpha) t} \|w(0)\|^{2} = e^{-2(\lambda + \alpha) t} \|v_2 - v_1\|^{2}.
    \end{align*}

\end{proof}

\section{Compact global attractors}\label{Sec4}

\begin{definition}\label{def2.7.3}
A cocycle $\varphi $ over $(Y,\mathbb T,\sigma)$ with the fibre
$W$ is said to be \emph{compactly dissipative} (respectively,
\emph{uniformly compact dissipative}) if there exits a nonempty
compact $K \subseteq W$ such that
\begin{equation}\label{eqGA2.7.8}
\lim_{t \to + \infty} \beta (U(t,y)M,K) =0
\end{equation}
for every $M \in C(W)$ and $y\in Y$ (respectively, uniformly with
respect to $y\in Y$), where $U(t,y):=\varphi(t,\cdot,y)$.
\end{definition}

Denote by
$$
\omega_{y}(K):=\bigcap\limits_{t\ge
0}\overline{\bigcup\limits_{\tau \ge
t}\varphi(\tau,K,\sigma(-\tau,y))},
$$
where $\varphi$ is a compactly dissipative cocycle and $K$ is a
compact subset appearing in (\ref{eqGA2.7.8}).

\begin{theorem}\label{thGA1}\cite{Che_1997},\ \cite[Ch.II]{Che_2015},\ \cite{Che_2022} Let $Y$ be a compact metric space
then the following statements are equivalent:
\begin{enumerate}
\item the cocycle $\langle W,\varphi,(Y,\mathbb T,\sigma)\rangle$
is uniformly compactly dissipative; \item the skew-product
dynamical system $(X,\mathbb T,\pi)$ ($X:=W\times Y, \pi
=(\varphi,\sigma)$) is compact dissipative.
\end{enumerate}
\end{theorem}

\begin{definition}\label{defCGA0_1} A family $I=\{I_{y}|\ y\in Y\}$ of
compact subsets $I_{y}$ of $W$ is said to be a compact global
attractor for the cocycle $\langle W,\varphi,(Y,\mathbb
R,\sigma)\rangle$ if the following conditions are fulfilled:
\begin{enumerate}
\item the set
\begin{equation}\label{eqCGA1}
\mathcal I :=\bigcup \{I_{y}|\ y\in Y\}\nonumber
\end{equation}
is precompact; \item the family of subsets $\{I_{y}|\ y\in Y\}$ is
invariant, i.e., $\varphi(t,I_{y},y)=I_{\sigma(t,y)}$ for all
$(t,y)\in \mathbb R_{+}\times Y$; \item
\begin{equation}\label{eqCGA2}
\lim\limits_{t\to +\infty}\sup\limits_{y\in
Y}\beta(\varphi(t,M,y),\mathcal I)=0\nonumber
\end{equation}
\end{enumerate}
for every compact subset $M$ from $W$.
\end{definition}

\begin{theorem}\label{thGA2} \cite{Che_2022} Let $Y$ be a compact metric space, $Y$ be invariant
(i.e., $\sigma(t,Y)=Y$ for all $t\in \mathbb T$) and $\varphi$ be
a cocycle over $(Y,\mathbb T,\sigma)$ with the fibre $W$. If the
cocycle $\varphi$ is uniformly compactly dissipative, then it has
a compact global attractor $I=\{I_{y}|\ y\in Y\}$, where
$I_{y}:=\omega_{y}(K)$ and the nonempty compact subset of $W$
appearing in the equality (\ref{eqGA2.7.8}).
\end{theorem}

\begin{definition}\label{defGA} Let $\langle W,\varphi,(Y,\mathbb R,\sigma)\rangle
$ be compactly dissipative, $K$ be the nonempty compact subset of
$W$ appearing in the equality (\ref{eqGA2.7.8}) and
$I_{y}:=\omega_{y}(K)$ for all $y\in Y$. The family of compact
subsets $\{I_y|\ y\in Y\}$ is said to be the \emph{Levinson
center} \cite[Ch.II]{Che_2020} (compact global attractor) of
nonautonomous (cocycle) dynamical system $\langle
W,\varphi,(Y,\mathbb R,\sigma)\rangle $.
\end{definition}

\begin{definition}\label{def3.11} Let
$\langle W,$ $\phi,$ $ (Y,$ $\mathbb R,$ $\sigma)\rangle $
(respectively, $(X,$ $\mathbb R_{+},$ $\pi))$ be a cocycle
(respectively, one sided dynamical system). A continuous mapping
$\nu :\mathbb R \to W $ (respectively, $\gamma :\mathbb R \to X$)
is called an entire trajectory of the cocycle $\phi $
(respectively, of the dynamical system $(X,\mathbb R_{+},\pi)$)
passing through the point $(u,y)\in W\times Y$ (respectively,
$x\in X$) for $t=0$ if $ \phi (t,\nu (s),\sigma(s,y))=\nu (t+s) \
\mbox{and} \ \nu (0)=u $ (respectively, $\pi (t,\gamma (s))=\gamma
(t+s)\ \mbox{and} \ \gamma(0)=x$) for all $ t\in \mathbb R_{+}$
and $s\in \mathbb R.$
\end{definition}

\begin{theorem}\label{thGAC1} \cite{Che_2022}, \cite[Ch.II]{Che_2024}
Let $\langle W,\varphi,(Y,\mathbb R,\sigma)\rangle $ be compactly
dissipative nonautonomous dynamical system, $\{I_{y}|\ y\in Y\}$
be its Levinson center $w\in I_{y}\ (y\in Y)$ if and only if there
exits a whole trajectory $\nu :\mathbb{R}\to W$ of the cocycle
$\varphi$ satisfying the following conditions: $\nu(0) = w$ and
$\nu(\mathbb R)$ is relatively compact.
\end{theorem}

\begin{definition}\label{defCGA1} A cocycle $\varphi$ is said to be
dissipative if there exists a bounded subset $K\subset \mathfrak
B$ such that for every bounded subset $B\subset \ \mathfrak B$ there
exists a positive number $L=L(B)$ such that
$\varphi(t,B,Y)\subseteq K$ for all $t\ge L(B)$, where
$\varphi(t,B,Y):=\{\varphi(t,u,y)|\ (u,y)\in B\times Y\}$.
\end{definition}

\begin{theorem}\label{thCGA1} \cite[Ch.II]{Che_2024} Assume that the metric space $Y$ is
compact and the cocycle $\langle \mathfrak B,\varphi,(Y,\mathbb
R,\sigma)\rangle$ is dissipative and asymptotically compact.
Then the cocycle $\varphi$ has a compact global attractor.
\end{theorem}

\begin{theorem}\label{thAC1} \cite{CS_2025} Assume that the following conditions hold:
\begin{enumerate}
\item[(D1)] the function $f\in C(\mathbb R,\ell_{2})$ is Lagrange
stable, i.e., $H(f)$ is a compact subset of $C(\mathbb
R,\ell_{2})$; \item[(D2)] the function $F\in (\mathbb R,\mathbb
R)$ is continuous differentiable and $F(0)=0$; \item[(D3)] there
are positive constants $\alpha$ and $\beta$ such that $sF(s)\le
-\alpha s^{2}+\beta$ for all $s\in \mathbb R$.
\end{enumerate}

Then the cocycle $\langle \ell_{2},\varphi,(H(f),\mathbb
R,\sigma)\rangle$ generated by the equation (\ref{eq2.3}) is
asymptotically compact.
\end{theorem}

\begin{theorem}\label{thCGA2} Under the Conditions
(\textbf{C1})-(\textbf{C3}) the equation (\ref{eq2.3}) (the
cocycle $\varphi$ generated by the equation (\ref{eq2.3})) has a
compact global attractor $\{I_{g}|\ g\in H(f)\}$.
\end{theorem}
\begin{proof}
The base space is $Y = H(f)$ is compact in \(C(\mathbb{R},
\ell_2)\) because the function $f$ is almost periodic. Consider
the equation:
\begin{equation}\label{eqCGA2_1}
u' = \nu \Lambda u - \lambda u + \widetilde{F}(u) + g(t), \quad g \in H(f).
\end{equation}
Note that
\begin{equation*}
    \frac{d}{dt} \frac{1}{2} \|u(t)\|^2 = \langle u(t), u'(t) \rangle = \langle u(t), \nu \Lambda u(t) -
     \lambda u(t) + \widetilde{F}(u(t)) + g(t) \rangle.
\end{equation*}
Split the inner product:
\[
\langle u(t), u'(t)\rangle = \nu \langle u(t), \Lambda u(t)\rangle
- \lambda \langle u(t), u(t)\rangle + \langle u(t),
\widetilde{F}(u(t))\rangle + \langle u(t), g(t)\rangle.
\]
Note that
\begin{itemize}
    \item[-]  \(\nu \langle u, \Lambda u \rangle = -\nu \|D^+ u\|^2 \le 0\), since \(\nu > 0\).
    \item[-]  \(-\lambda \langle u, u \rangle = -\lambda \|u\|^2\), with \(\lambda > 0\).
    \item[-]  by Lemma \ref{lAP1}, \(\langle u, \widetilde{F}(u) \rangle \le -\alpha \|u\|^2\) for all $u\in \ell_{2}$.
    \item[-]  \(\langle u, g(t) \rangle \le \|u\| \|g(t)\|\), and since \(H(f)\) is compact, \(\sup\limits_{g \in H(f)} \sup\limits_{t\in \mathbb R} \|g(t)\| \le M < \infty\).
\end{itemize}

Thus we obtain
\[
\frac{d}{dt} \frac{1}{2} \|u(t)\|^2 \le -(\lambda + \alpha)
\|u(t)\|^2 + M \|u(t)\|.
\]
Let \(z(t) = \|u(t)\|^2\), then
\[
\frac{d}{dt} z(t) \le -2(\lambda + \alpha) z(t) + M \sqrt{z(t)}.
\]
and, consequently,
\[
z(t) \le z(0) e^{-2(\lambda + \alpha)t} + \frac{M}{\lambda + \alpha} (1 - e^{-2(\lambda + \alpha)t}).
\]
for all $t\in \mathbb R_{+}$.

For every $\varepsilon >0$ there exists $T(\varepsilon,z(0))>0$
such that \(z(t) \le \frac{M}{\lambda + \alpha} + \varepsilon\)
for all \(t \ge T(\varepsilon,z(0))\).

Define
$$
K = \{ u \in \ell_2 \mid \|u\| \le \frac{M}{\lambda + \alpha} +
\varepsilon \}.
$$
For every bounded \(B\), choose \(L(B)\) such that \(\varphi(t, B,
H(f)) \subseteq K\) for \(t \ge L(B)\). Thus, the cocycle is
dissipative.

By Lemma \ref{lAP1} under the Conditions
(\textbf{C1})-(\textbf{C3}) the conditions (D1)-(D3) of Theorem
\ref{thAC1} are fulfilled and, consequently, the cocycle
\(\varphi\) is asymptotically compact. Thus the cocycle $\varphi$
is both dissipative and asymptotically compact, i.e, all
hypotheses of Theorem \ref{thCGA1} are satisfied. Consequently,
the cocycle \(\varphi\) generated by the equation (\ref{eq2.3})
admits a compact global attractor \(\mathcal{I} = \{I_g \mid g \in
H(f)\}\). This completes the proof.
\end{proof}

\section{Invariant sections of monotone nonautonomous lattice dynamical
systems}\label{Sec5}

Below we prove that under some conditions a nonautonomous
dynamical system generated by the equation (\ref{eqCGA2_1}) admits
an invariant continuous section.

Let $(Y,\mathbb{R},\sigma)$ be a two-sided dynamical system,
$(X,\mathbb{R_{+}}, \pi)$ be a semi-group dynamical system and
$h:X\to Y$ be a homomorphism of $(X,\mathbb{R_{+}},\pi)$ onto
$(Y,\mathbb{R},\sigma)$.

Let $\langle(X,\mathbb{R_{+}},\pi),(Y,\mathbb{R},\sigma),h\rangle$
be a nonautonomous dynamical system.

\begin{definition} A mapping $\gamma :Y\mapsto X$ is called a
continuous invariant section of nonautonomous dynamical system
$\langle (X,\mathbb R_{+},\pi),(Y,\mathbb R,\sigma),h\rangle$ if
the following conditions are fulfilled:
\begin{enumerate}
\item $h(\gamma(y))=y$ for all $y\in Y$; \item
$\gamma(\sigma(t,y))=\pi(t,\gamma(y))$ for all $y\in Y$ and $t\in
\mathbb R_{+}$; \item $\gamma$ is continuous.
\end{enumerate}
\end{definition}

\begin{remark}\label{remIS1} A continuous mapping $\gamma :Y\to X$
is an invariant section of the nonautonomous dynamical system
$\langle(X,\mathbb{R_{+}},\pi),(Y,\mathbb{R},\sigma),h\rangle$ if
and only if the set $\gamma(Y)$ is an invariant subset of
$(X,\mathbb R_{+},\pi)$.
\end{remark}

\begin{theorem}\label{thAP1}
\cite[Ch.I]{Che_2024},\cite[Ch.I]{Sch85} Suppose that $h$ is a
homomorphism of the dynamical system $(Y,\mathbb R,\sigma)$ into
$(X,\mathbb R_{+},\pi)$. If $y\in Y$ is a stationary
(respectively, $\tau$-periodic, quasi-periodic with the frequency
$\mu_{1}\ldots,\mu_{m}$ or almost periodic), then the point
$x:=h(y)$ is also stationary (respectively, $\tau$-periodic,
quasi-periodic with the frequency $\mu_{1},\ldots,\mu_{m}$ or
almost periodic).
\end{theorem}

\begin{lemma}\label{lCIS1} Let \(Y\) be a compact metric space and let $(X,\mathbb R_{+},\pi)$ and
$(Y,\mathbb R,\sigma)$ be two dynamical systems and $\gamma
:Y\mapsto X$ be a continuous invariant section. If the point $y\in
Y$ is stationary (respectively, $\tau$-periodic, quasi-periodic
with the frequency $\mu_1,\ldots,\mu_{m}$ or almost periodic),
then the point $x=\gamma(y)$ is also stationary (respectively,
$\tau$-periodic, quasi-periodic with the frequency
$\mu_1,\ldots,\mu_{m}$ or almost periodic).
\end{lemma}
\begin{proof} Let $\gamma: Y\to X$ be a continuous invariant section of the nonautonomous dynamical system
$\langle (X,\mathbb T_{1},\pi), (Y,\mathbb T_{2},\sigma),
h\rangle$ then $\gamma$ is a homomorphism of the dynamical system
$(Y,\mathbb T_{2},\sigma)$ into $(X,\mathbb T_{1},\pi)$. By
Theorem \ref{thAP1} if the point $y\in Y$ is stationary
(respectively, $\tau$-periodic, quasi-periodic with the frequency
$\mu_1,\ldots,\mu_{m}$ or almost periodic), then the point
$x=\gamma(y)$ is so.
\end{proof}

\begin{definition}\label{defFT1} A continuous function $\gamma :\mathbb R\to
X$(respectively, $\nu :\mathbb R\to W$) is said to be a full
trajectory of the semi-group dynamical system $(X,\mathbb
R_{+},\pi)$ (respectively, of the cocycle $\langle
W,\varphi,(Y,\mathbb R,\sigma)\rangle$) if
$\pi(t,\gamma(s))=\gamma(t+s)$ (respectively, if
$\varphi(t,\nu(s),\sigma(s,y))=\nu(t+s)$) for all $(t,s)\in
\mathbb R_{+}\times \mathbb R.$
\end{definition}

\begin{lemma}\label{lH1} Let $\langle W,\varphi,(Y,\mathbb
R,\sigma)\rangle$ be a cocycle over $(Y,\mathbb R,\sigma)$ with
the fibre $W$. Assume that the following conditions are fulfilled:
\begin{enumerate}
\item the point $y\in Y$ is stationary (respectively,
$\tau$-periodic, quasi periodic with the frequency
$\mu_1,\ldots,\mu_{m}$ or almost periodic); \item
$h:=pr_{2}:X:=W\times Y\to Y$ and $\gamma =(\nu,Id_{Y})$) is a
homomorphism of $(Y,\mathbb R,\sigma)$ into $(X,\mathbb
R_{+},\pi)$ (respectively, into $\langle W,\varphi,(Y,\mathbb
R,\sigma)\rangle$).
\end{enumerate}

Then the point $x=\gamma(y)$ (respectively, $v=\nu(y)$ and
$x=(v,y)$) is also stationary (respectively, $\tau$-periodic,
quasi periodic with the frequency $\mu_1,\ldots,\mu_{m}$ or almost
periodic).
\end{lemma}
\begin{proof} This statement directly follows from Lemma
\ref{lCIS1} because the mapping $\gamma =(\nu,Id_{Y})$ is a
homomorphism from $(Y,\mathbb R.\sigma)$ into skew-product
dynamical system $(X,\mathbb R_{+},\pi)$ ($X=W\times Y$ and $\pi
=(\varphi,\sigma)$).
\end{proof}

\begin{remark}\label{r3.1}
1. A continuous section $\gamma\in\Gamma(Y,X)$ is invariant if and
only if $\gamma\in\Gamma(Y,X)$ is a stationary point of the
semigroup $\{S^t\ \vert \ t\in \mathbb{R}_{+} \}$, where $S^t:
\Gamma(Y,X)\to \Gamma(Y,X)$ is defined by the equality
\begin{equation}\label{eqG1_1}
(S^t\gamma)(y):=\pi(t,\gamma(\sigma(-t,y)))
\end{equation}
for all $y\in Y$ and $t\in \mathbb R_+$.

2. Let $\langle W,\varphi, (Y,\mathbb R,\sigma)\rangle$ be a
cocycle, $(X,\mathbb R_{+},\pi)$ be the skew-product dynamical
system associated by cocycle $\varphi$ (i.e., $X:=W\times Y$ and
$\pi :=(\varphi,\sigma)$) and $\langle (X,\mathbb R_{+},\pi),(Y,$
$\mathbb R,$ $\sigma),h\rangle$ be the nonautonoous dynamical
system generated by $\varphi$ (i.e., $h:=pr_{2}:X\to Y$). Then the
following statements hold:
\begin{enumerate}
\item $\gamma \in \Gamma (Y,X)$ if and only if $\gamma
=(\nu,Id_{Y})$, where $\nu \in C(Y,W)$ and $Id_{Y}$ is the
identity mapping in $Y$; \item the equality (\ref{eqG1_1}) in this
case can be rewrite as follows
\begin{equation}\label{eqG1}
(S^t\gamma)(y)=(\Phi^{t}\nu(y),y)=(\varphi(t,\nu(\sigma(-t,y)),\sigma(-t,y)),y),\nonumber
\end{equation}
\begin{equation}\label{eqG1.1}
\Phi^{t}(\nu(y),y)=\varphi(t,\nu(\sigma(-t,y)),\sigma(-t,y))
\nonumber
\end{equation}
for all $y\in Y$.
\end{enumerate}
\end{remark}

\begin{theorem}\label{thAPM1} Let $\langle \ell_{2} ,\varphi, (H(f),\mathbb
R,\sigma)\rangle$ be the cocycle generated by the equation
(\ref{eq2.3}). Under the condition (\textbf{C1})-(\textbf{C3})
there exists a unique invariant section $\nu :H(f)\to \ell_{2}$ of
the cocycle $\varphi$ and
\begin{equation}\label{eqAPM_1}
\|\varphi(t,v,g)-\nu(\sigma(t,g))\|\le e^{-(\alpha
+\alpha)t}\|v-\nu(g)\| \nonumber
\end{equation}
for all $t\in \mathbb R_{+}$ and $v\in \ell_{2}$.
\end{theorem}
\begin{proof}
For every $t\in \mathbb R_{+}$ define a mapping $\Phi^{t}:
C(H(f),\ell_{2})\to C(H(f),\ell_{2})$ by the equality
\begin{equation}\label{eqAPM_11}
(\Phi^{t}\psi)(g):=\varphi(t,\psi(\sigma(-t,g)),\sigma(-t,g))\nonumber
\end{equation}
for every $g\in H(f)$.

Note that the family of mappings $\{\Phi^{t}\}$ acting on the
space $C(H(f),\ell_{2})$ possesses the following properties:
\begin{enumerate}
\item $\Phi^{0}=Id_{C(H(f),\ell_{2})}$; \item for every $t\in
\mathbb R_{+}$ the mapping $\Phi^{t}:C(H(f),\ell_{2})\to
C(H(f),\ell_{2})$ is continuous; \item $\Phi^{t}\circ
\Phi^{\tau}=\Phi^{t+\tau}$ for all $t,\tau\in \mathbb R_{+}$,
where by $\circ$ the compositions of two mappings is denoted.
\end{enumerate}

Since
\begin{eqnarray}\label{eqAPM_12}
&
\|\Phi^{t}\nu_{1}-\Phi^{t}\nu_{2}\|_{C(H(f),\ell_{2})}:=\nonumber
\\
& \max\limits_{g\in
H(f)}\|\varphi(t,\nu_{1}(\sigma(-t,g)),\sigma(-t,g))-\varphi(t,\nu_{2}(\sigma(-t,g)),\sigma(-t,g))\|=\nonumber
\\
& \max\limits_{g\in
H(f)}\|\varphi(t,\nu_{1}(g),g)-\varphi(t,\nu_{2}(g),g)\|\le \nonumber \\
& e^{-(\lambda +\alpha)t} \max\limits_{g\in
H(f)}\|\nu_{1}(g)-\nu_{2}(g)\|=e^{-(\lambda +\alpha)
t}\|\nu_{1}-\nu_{2}\|_{C(H(f),\ell_{2})}\nonumber
\end{eqnarray}
for all $t\in \mathbb R_{+}$ and $\nu_{1},\nu_{2}\in
C(H(f),\ell_{2})$ then the mapping $\Phi^{t}$ (for $t>0$) is a
contraction. Thus we have a commutative semigroup of continuous
mappings acting on the space $C(H(f),\ell_{2})$ and the mapping
$\Phi^{t_0}$ ($t_0>0$) is a contraction and, consequently (see,
for example, \cite[Ch.I]{DK_1970}), there exists a unique common
fixed point $\nu\in C(H(f),\ell_{2})$ of the semigroup
$\{\Phi^{t}\}_{t\ge 0}$. By Remark \ref{r3.1} (item 2(i)) the
mapping $\nu$ is a continuous invariant section of the cocycle
$\varphi$.

Finally, by Theorem \ref{thAP2} we have
\begin{equation}\label{eqAPM_13}
\|\varphi(t,v,g)-\nu(\sigma(t,g))\|=\|\varphi(t,v,g)-\varphi(t,\nu(g),g)\|\le
e^{-(\lambda +\alpha)t}\|v-\nu(g)\|\nonumber
\end{equation}
for all $(t,v,g)\in \mathbb R_{+}\times \ell_{2}\times H(f)$.
Theorem is completely proved.
\end{proof}

\begin{coro}\label{corAPM1} Under the Conditions (\textbf{C1})-(\textbf{C3}) the
equation (\ref{eq2.3}) has a unique almost periodic solution.
\end{coro}

\section{Convergent nonautonomous lattice dynamical
systems}\label{Sec6}

Let $\langle W,\varphi,(Y,\mathbb R,\sigma)\rangle$ (or shortly
$\varphi$) be a cocycle over dynamical system $(Y,\mathbb
R,\sigma)$ with the fibre $W$.

\begin{definition}\label{defCS1} A cocycle $\langle W,\varphi,(Y,\mathbb
R,\sigma)\rangle$ with the compact base space $Y$ is said to be
convergent if the following conditions are fulfilled:
\begin{enumerate}
\item the cocycle $\varphi$ admits a compact global attractor $I
=\{I_{y}|\ y\in Y\}$; \item for all $y\in Y$ the set $I_{y}$
consists of a single point $\{w_{y}\}$, i.e., $I_{y}=\{w_y\}$.
\end{enumerate}
\end{definition}

\begin{theorem}\label{thCGA_2} Under the Conditions
(\textbf{C1})-(\textbf{C3}) the equation (\ref{eq2.3}) (the
cocycle $\varphi$ generated by the equation (\ref{eq2.3})) is
convergent, i.e., it has a compact global attractor $ I =\{I_{g}|\
g\in H(f)\}$ such that for every $g\in H(f)$ the set $I_{g}$
consists of a single point.
\end{theorem}
\begin{proof}
Under the Conditions (\textbf{C1})-(\textbf{C3}) by Theorem
\ref{thCGA2} the equation (\ref{eq2.3}) admits a compact global
attractor $I =\{I_{g}|\ g\in H(f)\}$.

To finish the proof of Theorem we need to show that for every
$g\in H(f)$ the set $I_{g}$ consists of a single point. To this
end we note that by Theorem \ref{thAPM1} there exists an invariant
section $\nu$ of the cocycle $\varphi$ generated by
(\ref{eq2.3})). Since $\{I_{g}|\ g\in H(f)\}$ is the maximal
compact invariant set of the cocycle $\varphi$
\cite[Ch.II]{Che_2024}, $I'=\{I'_{g}|\ g\in H(f)\}$, where
$I'_{g}=\{\nu(g)\}$, is a compact invariant set of $\varphi$ and,
consequently, we have $\nu(g)\in I_{g}$ for every $g\in H(f)$. We
will show that $I_{g}=\{\nu(g)\}$ for all $g\in H(f)$. Assume that
it is not true then there exists $g_0\in H(f)$ such that the set
$I_{g_0}$ contains at least two different points $v_{i}$
($i=1,2$), Since the set $ I =\{I_{g}|\ g\in H(f)\}$ is invariant
then there exists two entire trajectories $\gamma_{i}$ ($i=1,2$)
such that $\gamma_{i}(0)=(v_{i},g)$ and
$\gamma_{i}(t)=(\nu_{i}(t),\sigma(t,g))\in J$ for all $t\in
\mathbb R$ (or equivalently, $\nu_{i}(0)=v_{i}$ and $\nu_{i}(t)\in
I_{\sigma(t,g)}$ for all $t\in \mathbb R$) and $i=1,2$. Since
$v_1\not= v_2$ and $\bigcup \{I_{g}|\ g\in H(f)\}$ is a precompact
set then we have
\begin{equation}\label{eqC1}
0<C:=\sup\limits_{t\in \mathbb R}|\nu_{1}(t)-\nu_{2}(t)|< +\infty
.
\end{equation}
Note that $(v_{i},y)=\gamma_{i}(0)=\pi(t,\gamma_{i}(-t))$ for all
$t\in \mathbb R_{+}$ and, consequently, we receive
\begin{equation}\label{eqC2}
v_{i}=\varphi(t,\nu_{i}(-t),\sigma(-t,g))
\end{equation}
for all $t\in \mathbb R_{+}$. From (\ref{eqC1}) and (\ref{eqC2})
we obtain
\begin{eqnarray}\label{eqC3}
&
\|v_1-v_2\|=\|\varphi(t,\nu_{1}(-t),\sigma(-t,g))-\varphi(t,\nu_{2}(-t),\sigma(-t,g))\|\le
\nonumber \\
& e^{-(\lambda +\alpha)t}\|\nu_{1}(-t)-\nu_{2}(-t)\|\le
e^{-(\lambda +\alpha)t}C
\end{eqnarray}
for all $t\in \mathbb R_{+}$. Passing to the limit in (\ref{eqC3})
as $t\to +\infty$ and taking into account (\ref{eqC1}) we obtain
$v_1=v_2$. The last equality contradicts to our assumption. The
obtained contradiction proves our statement.

Thus every set $I_{g}$ consists of a single point. Since
$\nu(g)\in I_{g}$ for all $g\in H(f)$ then $I_{g}=\{\nu(g)\}$.
Theorem is completely proved.
\end{proof}

\section{Applications}\label{Sec7}

Finally, we will give an example which illustrate our general
results.

\begin{example}\label{exGA1} {\rm Let $\{\omega_{i}\}_{i\in \mathbb Z}$ be a
sequence of the real numbers ($\omega_{i}\not= 0$ for all $i\in
\mathbb Z$). For every $i\in \mathbb Z$ we define a function
$f_{i}\in C(\mathbb R,\mathbb R)$ by the equality
\begin{equation}\label{eqGA1}
f_{i}(t):=\frac{\sin(\omega_{i}t)}{2^{|i|}} \nonumber
\end{equation}
for all $t\in \mathbb R$.

Note that the functions $f_{i}$ ($i\in \mathbb Z$) possess the
following properties:
\begin{enumerate}
\item
\begin{equation}\label{eqGA_0}
|f_{i}(t)|\le \frac{1}{2^{|i|}}
\end{equation}
for all $t\in \mathbb R$ and $i\in \mathbb Z$; \item
\begin{equation}\label{eqGA_1}
|f_{i}'(t)|\le \frac{|\omega_{i}|}{2^{|i|}}
\end{equation}
for all $t\in \mathbb R$ and $i\in \mathbb Z$.
\end{enumerate}

\begin{lemma}\label{lGA_1} For every $i\in \mathbb Z$ the
function $f_{i}$ is bounded and uniformly continuous on $\mathbb
R$.
\end{lemma}
\begin{proof} This statement directly follows from (\ref{eqGA_0})
and (\ref{eqGA_1}).
\end{proof}

\begin{lemma}\label{lGA1} The following statements hold:
\begin{enumerate}
\item $f(t)\in \ell_{2}$ for all $t\in \mathbb R$, where
$f(t):=(f_{i}(t))_{i\in \mathbb Z}$; \item for every $\varepsilon
>0$ there exists a number $n(\varepsilon)\in \mathbb N$ such that
\begin{equation}\label{eqGA2}
\sum\limits_{|i|>
n(\varepsilon)}|f_{i}(t)|^{2}<\frac{\varepsilon^{2}}{8}
\end{equation}
for all $t\in \mathbb R$.
\end{enumerate}
\end{lemma}
\begin{proof} Since
\begin{equation}\label{eqGA3}
\|f(t)\|^{2}=\sum\limits_{i\in \mathbb
Z}|f_{i}(t)|^{2}=\sum\limits_{i\in \mathbb
Z}\frac{\sin^{2}(\omega_{i}t)}{2^{2|i|}}\le \sum\limits_{i\in
\mathbb Z}\frac{1}{4^{|i|}}=\frac{11}{3} \nonumber
\end{equation}
then $f(t)\in \ell_{2}$ for all $t\in \mathbb R$.

Note that for every $\varepsilon >0$ there exists a number
$n(\varepsilon)\in \mathbb N$ such that
\begin{equation}\label{eqGA4}
\sum\limits_{|i|>
n(\varepsilon)}\frac{1}{4^{|i|}}<\frac{\varepsilon^{2}}{8} .
\end{equation}

By (\ref{eqGA_0}) and (\ref{eqGA4}) we obtain
\begin{equation}\label{eqGA_05}
\sum\limits_{|i|\ge n(\varepsilon)}|f_{i}(t)|^{2}\le
\sum\limits_{|i|\ge
n(\varepsilon)}\frac{1}{4^{|i|}}<\frac{\varepsilon^{2}}{8}
\end{equation}
for all $t\in \mathbb R$.
\end{proof}

Consider the function $f:\mathbb R\to \ell_{2}$ defined by
\begin{equation}\label{eqF_01}
f(t):=(f_{i}(t))_{i\in \mathbb Z}
\end{equation}
for every $t\in \mathbb R$.

\begin{lemma}\label{lGA2} The following statements hold:
\begin{enumerate}
\item the function $f:\mathbb R\to \ell_{2}$ is uniformly
continuous on $\mathbb R$; \item $f(\mathbb R)$ is a precompact
subset of $\ell_{2}$.
\end{enumerate}
\end{lemma}
\begin{proof}
For every $\varepsilon >0$ we choose $n(\varepsilon)\in \mathbb N$
such that (\ref{eqGA_05}) holds. Since the functions $f_{i}$
($|i|\le n(\varepsilon)$) are uniformly continuous then for
$\varepsilon$ there exists a positive number $\delta
=\delta(\varepsilon)$ such that $|t_1-t_2|<\delta$ implies (see
Lemma \ref{lGA_1})
\begin{equation}\label{eqGA6}
\sum\limits_{|i|\le
n(\varepsilon)}|f_{i}(t_1)-f_{i}(t_2)|^{2}<\frac{\varepsilon^{2}}{2}
.
\end{equation}

Indeed. If we assume that it is not true, then there are
$\varepsilon_0>0$, $\delta_{n}\to 0$ ($\delta_{n}>0$),
$\{t_{i}^{k}\}\ (i=1,2)$ with $|t_{1}^{k}-t_{2}^{k}|<\delta_{k}$
and
\begin{equation}\label{eqGA6.1}
\sum\limits_{|i|\le
n(\varepsilon_{0})}|f_{i}(t_{1}^{k})-f_{i}(t_{2}^{k})|^{2}\ge
\frac{\varepsilon_{0}^{2}}{2}
\end{equation}
for every $k\in \mathbb N$.

Logically two cases are possible:

1.  The sequence $\{t_{1}^{k}\}$ is bounded. In this case without
loss of generality we can assume that the sequence $\{t_{1}^{k}\}$
is convergent. Denote by
\begin{equation}\label{eqGA6.2}
\bar{t}_{1}:=\lim\limits_{k\to \infty}t_{1}^{k}.\nonumber
\end{equation}
It is clear that in this case the sequence $\{t_{2}^{k}\}$ also
converges and
\begin{equation}\label{eqGA6.3}
\bar{t}_{2}:=\lim\limits_{k\to \infty}t_{2}^{k}=\lim\limits_{k\to
\infty}t_{1}^{k} =\bar{t}_{1}.
\end{equation}

Passing to the limit in (\ref{eqGA6.1}) as $k\to \infty$ and
taking into account (\ref{eqGA6.3}) we obtain $\varepsilon_{0}=0$.
The last relation contradicts to the choice of the number
$\varepsilon_{0}$. The obtained contradiction proves our statement
in this case.

2. The sequence $\{t_{1}^{k}\}$ contains an unbounded subsequence
$\{t_{1}^{k_{m}}\}$. In this case without loss of generality we
can assume that $|t_{1}^{k}|\to +\infty$ as $k\to \infty$. Since
the functions $f_{i}\ (|i|\le n(\varepsilon_{0}))$ are Lagrange
stable then we can assume that the sequences
\begin{equation}\label{eqGA6.4}
\{\sigma(t_{1}^{k},f_{i})\}=\{f_{i}^{t_{1}^{k}}\}\nonumber
\end{equation}
are convergent in the space $C(\mathbb R,\mathbb R)$. Denote by
\begin{equation}\label{eqGA6.5}
\widetilde{f}^{1}_{i}:=\lim\limits_{k\to
\infty}\sigma(t_{1}^{k},f_{i})\nonumber
\end{equation}
then the sequences
$\{\sigma(t_{2}^{k},f_{i})\}=\{f_{i}^{t_{2}^{k}}\}$ ($|i|\le
n(\varepsilon_{0})$) also converge in $C(\mathbb R,\mathbb R)$.
Since $h^{k}:=t^{k}_{2}-t^{k}_{1}\to 0$ as $k\to \infty$ we will
have
\begin{eqnarray}\label{eqGA6.6}
& \widetilde{f}_{2}^{i}:=\lim\limits_{k\to
\infty}\sigma(t_{2}^{k},f_{i})=\lim\limits_{k\to
\infty}\sigma(t_{1}^{k}+h^{k},f_{i})=\nonumber \\
& \lim\limits_{k\to
\infty}\sigma(h^{k},\sigma(t_{1}^{k},f_{i}))=\lim\limits_{k\to
\infty}\sigma(t_{1}^{k},f_{i})=\widetilde{f}_{1}^{i}
\end{eqnarray}
for every $|i|\le n(\varepsilon_{0})$. From (\ref{eqGA6.6}) we
obtain
\begin{equation}\label{eqGA6.7}
\widetilde{f}^{i}_{2}(0)=\lim\limits_{k\to
\infty}f^{i}(t^{k}_{2})=\lim\limits_{k\to
\infty}f^{i}(t^{k}_{1})=\widetilde{f}^{i}_{1}(0)
\end{equation}
for every $|i|\le n(\varepsilon_{0})$. Passing to the limit in
(\ref{eqGA6.1}) as $k\to \infty$ and taking into account
(\ref{eqGA6.7}) we obtain $\varepsilon_{0}=0$ which contradicts to
our assumption. The obtained contradiction completes the proof of
our statement.

On the other hand we have
\begin{eqnarray}\label{eqGA7}
\|f(t_1)-f(t_2)\|^{2}=\sum\limits_{|i|\le
n(\varepsilon)}|f_{i}(t_1)-f_{i}(t_2)|^{2}+ \sum\limits_{|i|>
n(\varepsilon)}|f_{i}(t_1)-f_{i}(t_2)|^{2}
\end{eqnarray}
for all $t_1,t_2\in \mathbb R$. From (\ref{eqGA2}), (\ref{eqGA6})
and (\ref{eqGA7}) we receive
$$
\|f(t_1)-f(t_2)\|^{2}=\sum\limits_{|i|\le
n(\varepsilon)}|f_{i}(t_1)-f_{i}(t_2)|^{2}+ \sum\limits_{|i|>
n(\varepsilon)}|f_{i}(t_1)-f_{i}(t_2)|^{2}\le
$$
$$
\sum\limits_{|i|< n(\varepsilon)}|f_{i}(t_1)-f_{i}(t_2)|^{2}+
\sum\limits_{|i|>
n(\varepsilon)}2(|f_{i}(t_1)|^{2}+|f_{i}(t_2)|^{2})<
$$
$$
\frac{\varepsilon^{2}}{2}+
2(\frac{\varepsilon^{2}}{8}+\frac{\varepsilon^{2}}{8})=\varepsilon^{2}
$$
and, consequently, $\|f(t_1)-f(t_2)\|<\varepsilon$ for all
$t_1,t_2\in \mathbb R$ with $|t_1-t_2|<\delta$.

Let now $v$ be an arbitrary element of the set $f(\mathbb R)$ then
there exists a number $s\in \mathbb R$ such that $v=f(s)$. By
Lemma \ref{lGA1} (item (ii)) for every $\varepsilon
>0$ there exists a number $n(\varepsilon)\in \mathbb N$ such that
\begin{equation}\label{eqGA9}
\sum\limits_{|i|> n(\varepsilon)}|v_i|^{2}= \sum\limits_{|i|>
n(\varepsilon)}|f_i(s)|^{2}<\frac{\varepsilon^{2}}{8} \nonumber
\end{equation}
and, consequently, by Theorem 5.25 \cite[Ch.V,p.167]{LS_1974} the
subset $f(\mathbb R)$ of $\ell_{2}$ is precompact.
\end{proof}

\begin{coro}\label{corGA1} The function $f$ is Lagrange stable,
i.e., the set $H(f)$ is a compact subset of $C(\mathbb
R,\ell_{2})$.
\end{coro}
\begin{proof} This statement follows from Lemmas \ref{lAPF02} and
\ref{lGA2}.
\end{proof}

\begin{lemma}\label{lGA3} The function $f\in C(\mathbb
R,\ell_{2})$ defined by (\ref{eqF_01}) is almost periodic.
\end{lemma}
\begin{proof}
For every $\varepsilon >0$ we choose $n(\varepsilon)\in \mathbb N$
such that
\begin{equation}\label{eqGA5}
\sum\limits_{|i|> n(\varepsilon)}|f_{i}(t)|^{2}\le \sum_{|i|>
n(\varepsilon)}\frac{1}{4^{|i|}}<\frac{\varepsilon^{2}}{8}
\end{equation}
holds for all $t\in \mathbb R$. The function $F\in C(\mathbb
R,\mathbb R^{2n(\varepsilon)+1})$ defined by
$F(t):=(f_{i}(t))_{|i|\le n(\varepsilon)}$ is almost periodic
because the function $f_{i}\in C(\mathbb R,\mathbb R)$ is $
\frac{2\pi}{\omega_{i}}$ periodic. By above for the positive
number $\frac{\varepsilon}{\sqrt{2(2n(\varepsilon)+1)}}$ there
exists a relatively dense subset $\mathcal F(\varepsilon)$ of
$\mathbb R$ such that for every $\tau \in \mathcal F(\varepsilon)$
we have
\begin{equation}\label{eqF2}
|f_{i}(t+\tau)-f_{i}(t)|<\frac{\varepsilon}{\sqrt{2(2n(\varepsilon)+1)}}
\end{equation}
for all $t\in \mathbb R$ and $|i|\le n(\varepsilon)$. Note that
$$
\|f(t+\tau)-f(t)\|^{2}=\sum\limits_{|i|\le
n(\varepsilon)}|f_{i}(t+\tau)-f_{i}(t)|^{2} + \sum\limits_{|i|>
n(\varepsilon)}|f_{i}(t+\tau)-f_{i}(t)|^{2}\le
$$
\begin{equation}\label{eqF3}
\sum\limits_{|i|\le n(\varepsilon)}|f_{i}(t+\tau)-f_{i}(t)|^{2}
+2\sum\limits_{|i|>
n(\varepsilon)}(|f_{i}(t+\tau)|^{2}+|f_{i}(t)|^{2})
\end{equation}
for all $t\in \mathbb R$ and $\tau \in \mathcal F(\varepsilon)$.
From (\ref{eqGA5}), (\ref{eqF2}) and (\ref{eqF3}) we obtain
\begin{equation}\label{eqF4}
\|f(t+\tau)-f(t)\|^{2}<\varepsilon^{2}/2+\varepsilon^{2}/2=\varepsilon^{2}\nonumber
\end{equation}
for all $t\in \mathbb R$ and, consequently, the function $f\in
C(\mathbb R,\ell_{2})$ is almost periodic. Lemma is proved.
\end{proof}

Consider the system of differential equations
\begin{equation}\label{eqAP_2}
u_{i}' = \nu (u_{i-1} - 2u_i + u_{i+1}) - \lambda u_{i} + F(u_i) +
\frac{\sin (\omega_{i} t )}{2^{|i|}} \quad (i \in \mathbb{Z}),
\end{equation}
where $F(u) = -u(1+u^{2})$ for all $u \in \mathbb{R}$.

Along with the system of equations (\ref{eqAP_2}), consider the
(equivalent) equation
\begin{equation}\label{eqAP03}
\xi' = \Lambda \xi - \lambda \xi + \tilde{F}(\xi) + f(t)
\end{equation}
in the space $\ell_{2}$.

Taking into account the results above it is easy to check that the
Conditions (\textbf{C1})-(\textbf{C3}) are fulfilled and, consequently, the equation
(\ref{eqAP03}) has a unique almost periodic solution $\varphi(t)$
which is globally asymptotically stable.}
\end{example}

\section{Some generalizations}\label{Sec8}

In this section, we outline how some of the results in this
article can be generalized to almost automorphic (respectively,
recurrent) lattice dynamical systems.

\begin{definition}\label{def14.3}\rm
A point $x \in X $ is called almost recurrent, if for every
$\varepsilon >0$ the set $\mathcal F(x,\varepsilon):=\{\tau \in
\mathbb T|\ \rho(\pi(\tau,x),x)<\varepsilon \}$ is relatively
dense.
\end{definition}

\begin{definition}\label{def4.4}\rm
If a point $x\in X$ is almost recurrent and its trajectory
$\Sigma_x$ is pre-compact, then $x$ is called (Birkhoff)
recurrent.
\end{definition}

Denote by $\mathfrak N_{x}:=\{\{t_n\}\subset \mathbb R|\
\pi(t_n,x)\to x\ \mbox{as}\ n\to \infty \}$.

\begin{definition}\rm
A point $x\in X$ is called {\em Levitan almost periodic}
\cite{Lev-Zhi} (see also \cite{Bro79,Che_2008,Lev_1953}), if there
exists a dynamical system $(Y,\mathbb T,\sigma)$ and a Bohr almost
periodic point $y\in Y$ such that $\mathfrak N_{y}\subseteq
\mathfrak N_{x}.$
\end{definition}

\begin{definition} \rm
A point $x\in X$ is called {\em almost automorphic} if it is st.
$L$ and Levitan almost periodic.
\end{definition}

\begin{remark}\label{remAA1} 1. Every almost periodic motion is
almost automorphic.

2. There exists almost automorphic motions which are not almost
periodic (see, for example, \cite{CL_2020}).
\end{remark}

\begin{theorem}\label{thAA1} \cite{Che_2020} Let $h:Y\to X$ be a homomorphism of
dynamical system $(Y,\mathbb T,\sigma)$ into $(X,\mathbb T,\pi)$.
If the point $y\in Y$ is almost recurrent (respectively,
recurrent, Levitan almost periodic or almost automorphic), then
the point $x:=h(y)$ is so.
\end{theorem}

\begin{coro}\label{corAA1} Let $\gamma :Y\to X$be a continuous
invariant section of nonautonomous dynamical system $\langle
(X,\mathbb T,\pi),(Y,\mathbb T,\sigma),h\rangle$. If the point
$y\in Y$ is almost recurrent (respectively, recurrent, Levitan
almost periodic or almost automorphic), then the point
$x:=\gamma(y)$ is so.
\end{coro}
\begin{proof} This statement directly follows from Theorem
\ref{thAA1} because every section $\gamma$ of the nonautonomous
dynamical system $\langle (X,\mathbb T,\pi),(Y,\mathbb
T,\sigma),h\rangle$ is a homomorphism of the dynamical system
$(Y,\mathbb T,\sigma)$ into $(X,\mathbb T,\pi)$.
\end{proof}

\begin{theorem}\label{thAA2} Under the Conditions (\textbf{C1})-(\textbf{C3})
if the function $f$ is almost automorphic (respectively,
resurrent), then the equation (\ref{eq2.3}) has a unique almost
automorphic (respectively, recurrent) solution which is globally
exponentially stable.
\end{theorem}
\begin{proof} This statement follows from Theorem \ref{thAPM1} and
Corollary \ref{corAA1}.
\end{proof}

\section{Funding}

This research was supported by the State Programs of the Republic
of Moldova "Remotely Almost Periodic Solutions of Differential
Equations (25.80012.5007.77SE)" and partially was supported by the
Institutional Research Program 011303 "SATGED", Moldova State
University.

\section{Conflict of Interest}

The authors declare that they have not conflict of interest.

\medskip
\textbf{ORCID (D. Cheban):} https://orcid.org/0000-0002-2309-3823

\textbf{ORCID (A. Sultan):} https://orcid.org/0009-0003-9785-9291

\end{document}